\documentclass{amsart}
\usepackage{amscd,amsmath,amssymb,amsthm}
\newcommand{\N}{\mathbb{N}}
\newcommand{\Z}{\mathbb{Z}}

\newcommand{\R}{\mathbb{R}}
\newcommand{\E}{\mathbb{E}}
\renewcommand{\H}{\mathbb{H}}
\newcommand{\inv}{^{-1}}
\newcommand{\ab}{_{ab}}
\renewcommand{\span}[1]{\left<#1\right>}
\newcommand{\spanst}[2]{\span{\,#1\mid#2\,}}
\newcommand{\set}[1]{\left\{#1\right\}}
\newcommand{\setst}[2]{\set{\,#1\mid#2\,}}
\newcommand{\abs}[1]{\left|#1\right|}
\newcommand{\Abs}[1]{\left\|#1\right\|}
\newcommand{\onto}{\twoheadrightarrow}

\newcommand{\func}[4][\to]{#2\colon#3#1#4}
\newcommand{\lp}[1]{\Z\left[#1,#1\inv\right]}
\newcommand{\norm}{\vartriangleleft}
\newcommand{\ind}[2]{\left[#1:#2\right]}
\DeclareMathOperator{\head}{head}
\DeclareMathOperator{\tail}{tail}
\DeclareMathOperator{\tr}{tr}
\newtheorem{prop}{Proposition}
\newtheorem{lem}[prop]{Lemma}
\newtheorem{cor}[prop]{Corollary}
\theoremstyle{definition}
\newtheorem{defin}{Definition}
\theoremstyle{remark}
\newtheorem*{notat}{Notation}
\newtheorem*{rk}{Remark}
\begin{document}

\title{Deep pockets in lattices and other groups}
\date{\today}
\author[A.D.~Warshall]{Andrew D. Warshall}
\thanks{We thank our advisor, Andrew Casson, John Mackay and Joshua
Zelinsky for their helpful comments.}
\address{Yale University\\Department of Mathematics\\P.O. Box
208283\\New Haven, CT 06520-8283\\USA}
\email{andrew.warshall@yale.edu}
\begin{abstract}
We show the nonexistence of deep pockets in a large class of groups,
extending a result of Bogopol'ski\u{i}. We then give examples of
important groups (namely lattices in Nil and Sol) which have deep
pockets.
\end{abstract}
\maketitle

\section{Introduction}
Let $G$ be an infinite group and $A$ a finite generating set for
$G$. We define the \emph{depth} (or more verbosely the \emph{dead-end
depth}) of an element $g\in G$ with respect to $A$ to be the distance
(in the word metric with respect to $A$) from $g$ to the complement of
the radius-$d_A(1,g)$ closed ball about the identity in $G$. In fact
this definition makes sense if $G$ is replaced with an arbitrary
pointed metric space; this will be used in Section~\ref{euc}. Depth
must be finite for all $g\in G$ since $G$ is infinite and finitely
generated. If the depth of $g$ is $>1$ then $g$ is called a \emph{dead
end}. (The idea is that a geodesic word representing $g$ cannot be
extended further if it is to remain a geodesic.) We define the depth
of $G$ with respect to $A$ to be the supremum of the depths of all
elements of $G$.

A natural question to ask is whether the depth of $G$ for some given
generating set $A$ is finite; if it is not, $G$ is said to have
\emph{deep pockets} with respect to $A$. One result in this direction
has been proven by Bogopol'ski\u{i} (see \cite{B}), who showed that
all hyperbolic groups have finite depth with respect to all generating
sets. We extend this result to a broader class of groups, namely those
with a regular language of geodesics. To be precise, we prove the
following

\begin{prop}\label{regbound}
Let $G$ be an infinite group with a regular language $L$ of geodesics
with respect to $A$ a finite generating set for $G$. Then there exists
a uniform bound for the depth of elements of $G$ with respect to $A$.
\end{prop}

\begin{defin}
A group is said to be \emph{weakly geodesically automatic} if it
admits an automatic structure made up entirely of geodesics.
\end{defin}

\begin{cor}
Any group $G$ weakly geodesically automatic with respect to a
generating set $A$ has a uniform bound on depth with respect to $A$.
\end{cor}

Two cases of this corollary are hyperbolic groups and finitely
generated abelian groups, each with respect to any set of generators;
for details see \cite{ECHLPT}. The first of these coincides with
Bogopol'ski\u{i}'s result.

The next proposition provides a slight extension of the finitely
generated abelian case.

\begin{prop}\label{eucd}
Any euclidean group has a uniform bound on depth with respect to any
set of generators.
\end{prop}

These propositions seem unsurprising, since previous constructions of
groups with deep pockets tended to involve wreath products or the
like; see \cite{CT}, \cite{CR} and \cite{RW}. However, we here give
two examples of comparatively well-known and well-behaved groups whose
depth is infinite with respect to reasonably standard generating sets.

\begin{prop}\label{heisd}
Let $H=\spanst{a,b}{[a,[a,b]],[b,[a,b]]}$ be the discrete Heisenberg
group. Then $H$ has infinite depth with respect to its two-element
generating set $\set{a,b}$.
\end{prop}

\begin{prop}\label{sold}
Let $R$ be a hyperbolic automorphism of $\Z^2$ and let
$G_R=\Z^2\rtimes_R\Z$. Let $\set{a,b}$ be the standard generating set
for $\Z^2$ and $\set{c}$ be that for $\Z$. Then $G_R$ has arbitrarily
deep dead ends with respect to $\set{a,b,c}$.
\end{prop}

It is interesting that $H$ and $G_R$ are lattices in two of Thurston's
eight three-dimensional model geometries (see \cite{T}), namely Nil
and Sol respectively. In contrast, it follows from
Proposition~\ref{regbound} (or \cite{B}) and Proposition~\ref{eucd}
that no lattice in $\H^3$, $S^2\times\R$, $\widehat{SL_2(\R)}$ or
$\E^3$ can have deep pockets. Since the question of admitting a
lattice with deep pockets is meaningless for $S^3$ since it is compact
(so that any lattice must be finite), there remains one model geometry
for which it is (to the best of this author's knowledge) open, namely
$\H^2\times\R$.

\section{Proof of Proposition~\ref{regbound}}
\begin{proof}
Since $G$ is infinite, every individual element has finite depth.  Let
$n$ be the number of states of a (deterministic) finite state
automaton $F$ accepting $L$. Then it will suffice to prove the
existence of a uniform bound on depth for elements at distance at
least $n$ from the identity in $G$ with respect to $A$. Let $g$ be
such an element and consider the sequence of states assumed by $F$ as
it reads the element $w$ of $A$ representing $g$. This sequence has
length at least $n+1$, so, by the pigeonhole principle, the last $n+1$
terms of this sequence must contain at least one repetition. We thus
express $w$ as $abc$, where $\abs{bc}\le n$, $\abs{b}>0$ and $F$ is in
the same state after reading $ab$ as after reading $a$. Then $F$ will
also accept $abbc$, since after reading $abb$ it will be in the same
state as after reading $ab$, so that reading $c$ will get it into an
accept state. Since all words accepted by $F$ are geodesic, $abbc$ is
thus a geodesic word of length greater than that of $w=abc$ (since
$\abs{b}>0$) and representing a group element within
$\abs{c}+\abs{bc}\le2\abs{bc}\le2n$ of $g$. Thus the depth of $g$ is
$\le2n$, so we are done.
\end{proof}

\section{Proof of Proposition~\ref{eucd}}\label{euc}
We need the following auxiliary result, which will also be used in
Section~\ref{sol}.

\begin{prop}\label{fuzz}
Let $f$ be a function from a metric space $A$ to $\Z$ and $n\in\Z$. Suppose
there exists $a\in A$ and $r\in\Z$ such that for all $a'\in B_r(a)$
$f(a')\le f(a)+n$. Then there exists some $a'\in A$ such that $f$
attains a maximum on $B_{r/n}(a')$ at $a'$.
\end{prop}

\begin{proof}
Define $g(x)=\max_{b\in B_a(x)}f(b)$ for $x>0$, $f(a)$ for
$x=0$. Clearly, $g$ is nondecreasing and, by hypothesis, $g(r)\le
f(a)+n=g(0)+n$. Thus there must exist some $x$ such that $g$ is
constant on the interval $(x,x+r/n]$, where $x+r/n\le r$. Then we may
just set $a'$ to be an element of $A$ at distance $\le x$ from $a$
such that $f(a')=g(x)$. By the triangle inequality and the definition
of $g$, $f$ is bounded on $B_{r/n}(a')$ by $g(x+r/n)=g(x)=f(a')$.
\end{proof}

From this we deduce

\begin{prop}\label{pocketfuzz}
Suppose a pointed metric space $(A,a_0,d)$ has deep pockets. Then if
$C$ is a real number and $d'$ is a integer-valued metric on $A$ such
that
\[
\abs{d'(a_1,a_2)-d(a_1,a_2)}<C
\]
for all $a_1$, $a_2\in A$ then the resulting pointed metric space
$(A,a_0,d')$ will still have deep pockets.
\end{prop}

\begin{proof}
Let $\func{f}{A}{\Z}$ be distance from the origin with respect to the
new metric, and $a$ have depth at least $r+C$ with respect to
$(A,a_0,d)$. Then Proposition~\ref{fuzz} gives a dead end with respect
to $(A,a_0,d')$ of depth at least $r/C$. Since $r$ was arbitrary, we
are done.
\end{proof}

With these preliminaries out of the way, we approach the heart of the
proof, starting with the following

\begin{defin}
We define a \emph{weighted} generating set for $\Z^n$ to be an ordered
pair $(A,\mu)$, where $A$ is a generating set for $\Z^n$ and
$\func{\mu}{A}{\N}$ is called the \emph{weight} function. We say that
the \emph{length} of the word $a_1a_2\ldots a_m$ with respect to the
weight function $\mu$ is $\sum_{i=1}^m\mu(a_i)$.
\end{defin}

It is clear that a weighted generating set defines a word metric just
as does an ordinary generating set.

\begin{prop}\label{bweight}
Every finite weighted generating set for $\Z^n$ has depth bounded for
all group elements.
\end{prop}

Let $(A,\mu)$ be a generating set for $\Z^n$. Let $M$ be the least
common multiple of the elements of the image of $\mu$ and let $A'$ be
the set consisting of $aM/\mu(a)$ for every $a\in A$. If we regard
$\Z^n$ as being embedded in the obvious way as a lattice in $\R^n$ and
let $B$ be the convex hull of $A'\cup{A'}\inv$ in $\R^n$, then $B$
will clearly be a convex polytope. In fact, its dimension will be $n$,
since otherwise $A$ would be contained in an $n-1$-dimensional
subspace, hence could not generate. We have the following

\begin{prop}\label{geo}
If $a_1$, $a_2$, \ldots, $a_m\in A'\cup{A'}\inv$ are the vertices of a
facet of $B$ and $i_1$, $i_2$, \ldots, $i_m$ are nonnegative integers
then $w=a_1^{i_1}a_2^{i_2}\ldots a_m^{i_m}$ is a geodesic word. In
fact, $w$ lies on a geodesic ray.
\end{prop}

\begin{proof}
The facet of $B$ of which $a_1$, $a_2$, \ldots, $a_m$ are vertices
must lie in some hyperplane $H$ in $\R^n$, say
$\mathbf{a}\cdot\mathbf{x}=b$. Since $B$ is symmetric about the origin
and lies in one of the closed half-spaces bounded by $H$, if the
origin lay in $H$ then $B$ would lie entirely in $H$, so so would
$A'\cup{A'}\inv$, so $A\cup A\inv$ would as well, which is a
contradiction since $A$ generates $\Z^n$. Thus $b\ne0$, so we may
choose $\mathbf{a}$ so that $b=1$, so $H$ is given by the equation
$\mathbf{a}\cdot\mathbf{x}=1$ and $B$ lies in the closed half-space
given by $\mathbf{a}\cdot\mathbf{x}\le1$. In particular, so does
$A'\cup{A'}\inv$, so if $g\in\Z^n$ is given by a word of length $l$
with respect to $(A,\mu)$ then $\mathbf{a}\cdot g\le l/M$. Setting
$g=a_1^{i_1}a_2^{i_2}\ldots a_m^{i_m}$, we have $\mathbf{a}\cdot
g=\sum_{j=1}^m i_j$, so the minimal length for a word representing $g$
is $M\sum_{j=1}^m i_j$. But this is the length of $w$ since each $a_i$
has length $M$, proving the first sentence. The second sentence
follows clearly.
\end{proof}

\begin{prop}
Every element of $\Z^n$ is within bounded distance (with respect to
$(A,\mu)$) of a geodesic ray.
\end{prop}

\begin{proof}
By choosing a finite triangulation of each facet $f_i$ of $B$, we can
assume WLOG that they are all $n-1$-simplices, hence have $n$
vertices. As in the proof of Proposition~\ref{geo}, no facet may be
confined to any hyperplane through the origin, so the vertices of each
facet (since they clearly cannot lie in one $n-2$-plane) must be
linearly independent, hence generate a maximal-rank lattice. If we
regard this lattice as a based lattice embedded in $\R^n$ with the
usual norm, it will partition $\R^n$, hence $\Z^n$, into closed
parallelepipeds. Clearly there exists some bound $d_i$ on the
$A$-distance between any element of $\Z^n$ contained in the closed
unit parallelepiped generated by the facet $f_i$ and the origin. By
translation, we find that every element of $\Z^n$ is within $d_i$ of
some corner of every parallelepiped in which it lies. Since $B$ has
finitely many facets, we may set $d$ to be the greatest of the $d_i$.

I claim that every element of $\Z^n-B$ is within $d$ of a geodesic
ray. So let $g\in\Z^n-B$, and consider the line segment (in $\R^n$)
$s$ between $g$ and the origin. Since the origin is in $B$, $g\notin
B$ and $B$ is convex, $s$ intersects the boundary of $B$ exactly once
(say at $t$), so we let $a_1$, $a_2$, \ldots, $a_n$ be the vertices of
a (closed) facet containing this intersection. By the last paragraph,
$g$ is within distance $d$ of some corner of every parallelepiped of
the $a_1a_2\ldots a_n$-lattice to which it belongs. But the
coordinates of $g$ with respect to $a_1$, $a_2$, \ldots, $a_n$ as an
element of $\R^n$ must be nonnegative since they are obtained by
multiplying a positive number by the coordinates of $t$, which is (as
a vector) a weighted average of $a_1$, $a_2$, \ldots, $a_n$. Hence $g$
belongs to a parallelepiped all of whose corners have all their
coordinates nonnegative, hence (by Proposition~\ref{geo}) lie on a
geodesic ray, proving our claim.

Since $B$ is bounded, the case $g\in B$ is trivial.
\end{proof}

\begin{proof}[Proof of Proposition~\ref{bweight}]
Let $D$ be the bounded distance whose existence was proven in the
preceding proposition. Let $g\in\Z^n$ with minimal word length
$l$. Then $g$ is at distance at most $D$ from an element on a geodesic
ray whose $(A,\mu)$-distance from the origin is, by the triangle
inequality, at least $l-D$. Then there exists a group element at
distance at least $D+1$ but at most $D+M+1$ further out along this ray
whose distance from the origin is at least $l+1$ and, by the triangle
inequality again, at distance at most $2D+M+1$ from $g$. Thus $g$ has
depth at most $2D+M+1$, as claimed.
\end{proof}

Let $A$ be a generating set on $E$ a euclidean group and let
$\Z^n\norm E$. Let $\func[\onto]{\phi}{F_A}{E}$ be the projection
map. We denote the length of a word $w\in F_A$ by $l(w)$. Let $A'$ be
the set of words $w$ in $A\cup A\inv$ such that $\phi(w)\in\Z^n$ but,
for any proper subword $w'$ of $w$, $\phi(w')\notin\Z^n$. Let
$(B,\mu)$ be the finite weighted generating set for $\Z^n$ where
$B=\setst{a^g}{a\in A'}{g\in E}=\bigcup_{h\in E/\Z^n}{A'}^h$ and $\mu$
is defined so that $\mu(a^g)$ is the length of $a$ (as a word in $A$)
for $a\in A'$ and $g\in E$. It will follow from
Proposition~\ref{eucupper} that this is in fact a generating set. In
any event, let $\func{\pi}{\Z^B}{\Z^n}$ be the projection and let
$l_\mu(v)$ denote the length of the word $v\in\Z^B$ with respect to
$\mu$. We denote the distance of $x\in\Z^n$ from the origin (at the
moment, possibly infinite) with respect to this generating set by
$\Abs{x}$.

\begin{prop}\label{eucupper}
Suppose $w\in F_A$ with $\phi(w)=x\in\Z^n$. Then $\Abs{x}\le l(w)$. In
particular, $B$ is a generating set for $\Z^n$.
\end{prop}

\begin{proof}
Since $\phi(w)\in\Z^n$, $w$ must contain some minimal subword under
inclusion in $\phi\inv(\Z^n)$, say $w_1$. By definition, $w_1\in
A'$. Then $w=w_lw_1w_r$ may be replaced with $w_1^{w_l\inv}w_lw_r$,
where $w_1^{w_l\inv}\in B$ and
\[
\mu(w_1^{w_l\inv})=\mu(w_1)=l(w_1)=l(w)-l(w_l)-l(w_r)\le
l(w)-l(w_lw_r)
\]
by definition and the triangle inequality. The first claim is now
proven by induction. The second claim follows easily.
\end{proof}

\begin{prop}\label{euclower}
There exists a $D$ such that for any $x\in\Z^n$ there is a word $w\in
F_A$ with $\phi(w)=x$ and $l(w)\le\Abs{x}+D$.
\end{prop}

\begin{proof}
Let $m=\ind{E}{\Z^n}$ and pick a set of words $w_1$, \dots, $w_m$ in
$F_A$ representing each of the finitely many cosets of $\Z^n$ in
$E$. Consider a word $v\in\Z^B$ with $\pi(v)=x$ and
$l_\mu(v)=\Abs{x}$. Since $\Z^n$ is abelian, we can assume that $v$ is
ordered to separate out the letters in each ${A'}^h$. The word $v$ is
thus divided into $m$ pieces $v_1v_2\dots v_m$, where each
$v_i\in\Z^{{A'}^h}$ for some specific $h\in E/\Z^n$. Thus each $v_i$
can be replaced with $w_i\inv v_i'w_i$ for some $v_i'\in\Z^{A'}$ with
$l_\mu(v_i')=l_\mu(v_i)$. But then each $v_i'$ can be replaced with a
$v_i''\in F_A$ with $l(v_i'')\le l_\mu(v_i')=l_\mu(v_i)$ by
construction of $A'$ and $\mu$. We thus have
\[
l(w_i\inv v_i''w_i)\le2l(w_i)+l(v_i'')\le2l(w_i)+l_\mu(v_i).
\]
Concatenating these yields $w=w_1\inv v_1''w_1w_2\inv v_2''w_2\dots
w_m\inv v_n''w_m\in F_A$ with $\phi(w)=x$ and
\begin{multline*}
l(w)\le\sum_{i=1}^ml(w_i\inv
v_i''w_i)\le\sum_{i=1}^m(2l(w_i)+l_\mu(v_i))\\=2\sum_{i=1}^ml(w_i)+l_\mu(v)=2\sum_{i=1}^ml(w_i)+\Abs{x},
\end{multline*}
so we are done if we set $D=2\sum_{i=1}^ml(w_i)$.
\end{proof}

\begin{proof}[Proof of Proposition~\ref{eucd}]
It follows from Propositions~\ref{bweight}, \ref{eucupper},
\ref{euclower} and \ref{pocketfuzz} that $\Z^n$ cannot have deep
pockets with respect to the subspace metric induced from the metric
with respect to $A$ on $E$. But since $\ind{E}{\Z^n}<\infty$, every
element of $E$ lies within bounded distance $b$ of $\Z^n$. Thus if $E$
had a dead end of depth $d$ then $\Z^n$ would have (with respect to
the subspace metric) a dead end of depth at least $d/b+1$, by
Proposition~\ref{fuzz}. Thus $E$ cannot have deep pockets, as claimed.
\end{proof}

\section{Proof of Proposition~\ref{heisd}}
We now turn our attention to the Heisenberg group $H$ with respect to
its standard two-generator presentation
$\spanst{a,b}{[a,[a,b]],[b,[a,b]]}$. (We use $[x,y]$ to denote $x\inv
y\inv xy$.) It is well known that every element of $H$ is expressible
uniquely as a word of the form $a^ib^j[a,b]^k$, $a$, $b$, $c\in\Z$. In
general, of course, these words are not of minimal length; for
example, $ba$ is expressed as $ab[a,b]\inv$. However, we give the
following way of visualizing elements of $H$, from which insight can
be gained as to minimal-length words.

Consider $i$ and $j$ as coordinates in the standard integer lattice
$\Z^2$. Then we can represent a word in $a$ and $b$ as a path in
$\Z^2$ starting at the origin. Multiplication on the right by
$a^{\pm1}$ corresponds to appending to the path a segment one unit to
the right or left, while multiplication by $b^{\pm1}$ corresponds to
appending a segment one unit up or down. It is clear that this gives a
one-to-one correspondence between words and paths and that a word
representing $a^ib^j[a,b]^k$ corresponds to a path from the origin to
$(i,j)$. If to this path we append a vertical segment leading from
$(i,j)$ to $(i,0)$ and a horizontal segment from $(i,0)$ to the
origin, we get a loop based at the origin. This loop will correspond
to a word representing
$a^ib^j[a,b]^kb^{-j}a^{-i}=a^ib^jb^{-j}a^{-i}[a,b]^k=[a,b]^k$. We have
thus reduced the study of the correspondence between paths starting at
the origin and elements of $H$ to that of the correspondence between
loops based at the origin and powers of $[a,b]$.

\begin{prop}\label{invar}
If two loops based at the origin correspond to the same power of
$[a,b]$, then they enclose the same oriented area.
\end{prop}

\begin{proof}
Consider two loops based at the origin that correspond to the same
power of $[a,b]$. Then one can be transformed into the other by
appending loops conjugate to $a\inv a$, $aa\inv$, $b\inv b$, $bb\inv$,
$[a,[a,b]]$, $[b,[a,b]]$, $[[a,b],a]$ and $[[a,b],b]$. But the area
enclosed by the concatenation of two loops is the sum of the areas
enclosed by each one, so it suffices to show that any loop conjugate
to one of the eight listed above encloses zero area. Since conjugation
does not affect the area enclosed by a loop, it in fact suffices to
show that the eight loops listed above enclose zero area. But this can
be checked trivially.
\end{proof}

\begin{prop}\label{loopk}
Any loop based at the origin corresponds to $[a,b]^k$, where $k$ is
the oriented area it encloses.
\end{prop}

\begin{proof}
By Proposition~\ref{invar} and the discussion preceding it, it
suffices to show that the path corresponding to the word $[a,b]$
encloses an oriented area of $1$. But this is clear.
\end{proof}

\begin{prop}\label{length}
Any path $P$ from the origin to $(i,j)$ corresponds to
$a^ib^j[a,b]^k$, where $k$ is the oriented area enclosed by the
concatenation of $P$ and the path corresponding to $b^{-j}a^{-i}$.
\end{prop}

\begin{proof}
This is just the combination of Proposition~\ref{loopk} and the
discussion preceding Proposition~\ref{invar}.
\end{proof}

\begin{prop}\label{dd}
The element $[a,b]^{n^2+1}$ is at distance $4n+2$ from the identity.
\end{prop}

\begin{proof}
Clearly, $[a,b]^{n^2+1}$ is within $4n+2$ of the identity, since it
equals 
\[
a^{-n-1}b\inv ab^{-n+1}a^nb^n.
\]

So consider a word $w$ representing $[a,b]^{n^2+1}$. By
Proposition~\ref{length}, we know that $w$ corresponds to a loop based
at the origin enclosing oriented area $n^2+1$. If we let $x$ be the
measure of the projection of the loop to the horizontal axis and $y$
the measure of that to the vertical axis, we have $xy>n^2$, so, by the
arithmetic mean-geometric mean inequality, $x+y>2n$, so the length of
$w$ is $>4n$. Since $w$ corresponds to a loop, its length must be
even, so in fact it is at least $4n+2$.
\end{proof}

\begin{prop}\label{nd}
The element $a^ib^j[a,b]^k$ with $\abs{i}<n$, $\abs{j}<n$,
$\abs{k}<n(n+1)$ is within $4n+2$ of the identity.
\end{prop}

\begin{proof}
Since, by Proposition~\ref{length}, rotations of the plane by $\pi/2$
about the origin and reflections about the line $j=0$ do not affect
distance from the identity, we may assume WLOG that $i\ge\abs{j}$ and
$k\ge0$. Let $k=q(n+1)+r$, $0\le r\le n$, $0\le q\le n-1$. Then
clearly $b^{-q-1}a^rba^{n+1-r}b^qa^{-n-1}$ is a loop based at the
origin of length $2q+2n+4\le4n+2$ bounding oriented area
$n(n+1)$. Thus $b^{-q-1}a^rba^{n+1-r}b^qa^{i-n-1}b^j$ represents
$a^ib^j[a,b]^k$, and it will have length $\le 4n+2$ since $\abs{j}\le
i\le n+1$.
\end{proof}

\begin{prop}\label{nh}
Every element within distance $m$ of $[a,b]^{n^2+1}$ is of form
\[
a^ib^j[a,b]^k
\]
 with $\abs{i}\le m$, $\abs{j}\le m$ and $\abs{k}\le n^2+1+m(m-1)/2$.
\end{prop}

\begin{proof}
We proceed by induction. The statement is clear for $m=0$. So consider
$w=a^ib^j[a,b]^k$ with $\abs{i}\le m$, $\abs{j}\le m$ and
$\abs{k}<n^2+1+m(m-1)/2$. Then $wb=a^ib^{j+1}[a,b]^k$,
$wb\inv=a^ib^{j-1}[a,b]^k$, $wa=a^{i+1}b^j[a,b]^{k-j}$ and
$wa\inv=a^{i-1}b^j[a,b]^{k+j}$, so, by the triangle inequality, these
four new words can all be expressed as $a^{i'}b^{j'}[a,b]^{k'}$ with
$\abs{i'}\le m+1$, $\abs{j'}\le m+1$ and $\abs{k'}\le
n^2+1+m(m-1)/2+\abs{j}\le n^2+1+m(m-1)/2+m=n^2+1+m(m+1)/2$, as
claimed.
\end{proof}

\begin{prop}
If $n>2$, the element $g_n=[a,b]^{n^2+1}$ is a dead end of depth at
least $\sqrt{2n-4}+1$. In particular, $H$ has arbitrarily deep dead
ends with respect to the generating set $\set{a,b}$.
\end{prop}

\begin{proof}
By Proposition~\ref{dd}, $g_n$ is at distance $4n+2$ from the
identity. But, by Propositions~\ref{nd} and \ref{nh}, any group
element within $\sqrt{2n-4}$ of $g_n$ is within $4n+2$ of the
identity. The second sentence then follows clearly.
\end{proof}

\section{Proof of Proposition~\ref{sold}}\label{sol}
Let $R$ be a hyperbolic automorphism of $\Z^2$. We turn our attention
to the group
\[
G_R=\Z^2\rtimes_R\Z.
\]
We denote the standard generators of $\Z^2$ by $a$ and $b$ and that of
$\Z$ by $c$. This group is soluble, as its commutator is contained in
the span of $a$ and $b$, which commute with each other. Thus
$G_R'=\span{a,b}\cong\Z^2$ is a normal subgroup and
$G_R/G_R'=\span{G_R'c}$. The powers of $c$ are thus a complete set of
coset representatives, so every element of $G_R$ can be expressed
uniquely in the form $a^ib^jc^k$, where $i$, $j$ and $k\in\Z$. We will
use $a$, $b$ and $c$ to denote these three elements of $G_R$
throughout.

Let $F_S$ be the free group on the set $S$ and let $N$ be the normal
closure of $a$ and $b$ in $F_S$. We know $N$ is freely generated by
$\setst{a^{c^i},b^{c^i}}{i\in\Z}$. Any element $g\in F_S$ can be
expressed uniquely as a product $uc^i$, where $u\in N$ and
$i\in\Z$. Define $\func{\phi}{F_S}{N\ab}$ to send $g$ to the image of
$u$ under the abelianization map. (We denote by $N\ab$ the
abelianization of the group $N$.) Note that $N\ab$ is freely generated
as an abelian group by $\setst{a^{c^i},b^{c^i}}{i\in\Z}$.

Define $\func{\pi}{N\ab}{G_R'}$ to send $v\in N\ab$ to the element of
$G_R'$ represented by $v$. Similarly, define
$\func{\sigma_R}{F_S}{G_R}$ to send $w\in F_S$ to the element of $G_R$
represented by $w$.

\begin{notat}
If $w\in F_S$ or $N\ab$, let $l(w)$ denote its length in the given
sets of generators.
\end{notat}

We can now state the following proposition, which is central to the
proof of Proposition~\ref{sold} and whose proof we postpone.

\begin{prop}\label{bdiff}
Let $G_R=\Z^2\rtimes_R\span{c}$, where $R$ is a hyperbolic
automorphism of $\Z^2$, and let $S=\set{a,b,c}$, where $a$ and $b$ are
the standard generators of $\Z^2$. Then there exists $D$ with the
following property. Let $g=uc^n\in G_R$, where $u\in G_R'$ and
$n\in\Z$. Let $v\in N\ab$ be a minimal-length element of
$\pi\inv(u)$. Then there is some $w\in F_S$ with $\phi(w)=v$,
$\sigma_R(w)=g$ and
\[
l(w)-D\le\abs{g}\le l(w),
\]
where by $\abs{g}$ we mean the length of the minimal-length element of
$\sigma_R\inv(g)$.
\end{prop}

\begin{defin}
Let $g=uc^z\in G_R$ with $u\in G_R'$. Let $\Abs{g}$ represent the
minimal length of any element $w$ of $\sigma_R\inv(g)$ such that
$\phi(w)$ is a minimal-length element of $\pi\inv(u)$.
\end{defin}

\begin{rk}
Any element $w\in\sigma_R\inv(g)$ must satisfy $\pi(\phi(w))=u$.
\end{rk}

In this notation, Proposition~\ref{bdiff} implies that there exists a
$D$ dependent only on $R$ such that, for any $g\in G_R$,
$\Abs{g}-D\le\abs{g}\le\Abs{g}$.

\begin{notat}
Let $u_1$ denote $a$ and $u_2$ denote $b$.
\end{notat}

The following lemma, an adaptation of Cleary and Taback's result on
word length in wreath products, is used both in proving and in
applying Proposition~\ref{bdiff}.

\begin{lem}\label{ll}
Let $g=uc^z\in G_R=\Z^2\rtimes_R\span{c}$, where $R$ is a hyperbolic
automorphism of $\Z^2$ and $u\in G_R'$. Let
$v=\sum_{j=1}^nu_{k_j}^{c^{i_j}}\in N\ab$ (where all the $k_j$ are $1$
or $2$) such that $\pi(v)=u$ and let $S=\set{a,b,c}$ where $a$ and $b$
are the standard generators of $\Z^2$. Then the minimal length of all
words $w\in F_S$ such that $\sigma_R(w)=g$, $\phi(w)=v$ is
\begin{multline*}
2(\max(i_j,0)-\min(i_j,0))\\+\min(\abs{z-\max(i_j,0)}-\max(i_j,0),\abs{z-\min(i_j,0)}+\min(i_j,0))+n.
\end{multline*}
\end{lem}

\begin{proof}
Consider the natural maps $\pi_S$ and $\pi_T$ mapping $F_S$ and
$\Z^T$, respectively, to
\[
\Z^2\wr\Z=\spanst{a,b,c}{[a^{c^i},b^{c^j}],[a^{c^i},a^{c^j}],[b^{c^i},b^{c^j}]}.
\]
A word $w\in F_S$ satisfies $\sigma_R{w}=g$ and $\phi{w}=v$ if and
only if $\pi_S(w)=\pi_T(v)c^z$ in $\Z^2\wr\Z$. (This is because
$\phi(w)=v$ if and only if $\pi_S(w)=\pi_T(v)c^k$ for some $k$, and
then we have $\sigma_R(w)=\pi(v)c^k$.) Thus the minimal length for $w$
is simply the minimal length of a word in $F_S$ representing
$\pi_T(v)c^z$ in $\Z^2\wr\Z$. The result now follows from Cleary and
Taback's length formula in \cite{CT} for wreath products $G\wr\Z$, $G$
some group, with respect to a generating set consisting of generators
of $G$ and a standard generator for $\Z$.
\end{proof}

The final main ingredient is the following proposition, which tells us
about the geometry of the $\Abs{\cdot}$ norm and whose proof we
likewise postpone. Note that, for any group $G$ with generating set
$A$, $d_A$ refers to the word metric on $G$ induced by $A$.

\begin{prop}\label{flat}
Let $G_R=\Z^2\rtimes\Z$, where $R$ is a hyperbolic automorphism of
$\Z^2$, and let $a$ and $b$ be the standard generators of $\Z^2$. Let
$n$ be a positive integer. Then for all sufficiently large positive
integers $m$ there exists $B_{m,n}\subset G_R'$ such that every $u\in
B_{m,n}$ has the following properties:

\begin{itemize}
\item If $u$ is expressed as $a^ib^j$, $i$, $j\in\Z$, then $\abs{i}$,
$\abs{j}<\abs{\tr R}^m$ and
\item If $v\in N\ab$ is a minimal-length element of $\pi\inv(u)$ then
$M(v)>n$ and $m(v)<-n$.
\end{itemize}

Furthermore, there exists $p=p(R,m,n)\in\Z$ and $u\in B_{m,n}$ such
that for all $u'$ with $d_{\set{a,b}}(u.u')\le p$ the following
conditions hold:

\begin{itemize}
\item $g'\in B_{m,n}$ and
\item $\Abs{u'}\le\Abs{u}$.
\end{itemize}

Finally, $\lim_{m\to\infty}p(R,m,n)=\infty$ for every $R$ and $n$.
\end{prop}

With these ingredients, we give the

\begin{proof}[Proof of Proposition~\ref{sold}]
By Proposition~\ref{bdiff}, there exists $D$ such that, for any $g\in
G_R$, $\Abs{g}-D\le\abs{g}\le\Abs{g}$. Thus to prove that $G_R$ has
deep pockets with respect to the word metric in $S$ it suffices to
prove that for any $r$ we can find an element $h\in G_R$ such that
$\Abs{}$ attains a maximum on $B_r(h)$ at $h$, for we may then apply
Proposition~\ref{fuzz} with $f(g)=\abs{g}$ and $n=D$.

To that end, we choose $n\ge r$ and then $m$ such that $r\tau^r\le
p(R,m,n)$, where $p$ is defined as in Proposition~\ref{flat}; this
choice of $m$ is possible since
$\lim_{m\to\infty}p(R,m,n)=\infty$. Define $u\in G_R'$ to be the
element given by Proposition~\ref{flat} with $B_{m,n}$ as input. Then
let $g=u'c^i\in G_R$, $u'\in G_R$, be any element within distance $r$
of $u$ with respect to $S$. Then $\abs{i}\le r$. Furthermore, if $s\in
G_R'$, $sc^{i'}a^{\pm1}=sR^{-i'}(a^{\pm1})c^{i'}$ and similarly if $a$
is replaced by $b$, so right multiplication of any $u''c^{i'}$,
$u''\in G_R'$, $\abs{i'}\le r$, by any letter of $S$ or its inverse
changes $u''$ (with respect to $\set{a,b}$) by at most
$\abs{\tau}^r$. Thus $d_{\set{a,b}}(u,u')<r\tau^r\le p(R,m,n)$, so
$u'\in B_{m,n}$. By Proposition~\ref{flat} and Lemma~\ref{ll}, this
implies that $\Abs{g}\le\Abs{u'}$. But Proposition~\ref{flat} also
tells us that this is $\Abs{u'}\le\Abs{u}$, so we are done.
\end{proof}

It remains to prove Propositions~\ref{bdiff} and \ref{flat}. Let $p_R$
be the characteristic polynomial of $R$. We define
$\Lambda_R=\Z[t]/(p_R(t))$. There is a natural map
$\func[\onto]{\alpha_R}{\lp{t}}{\Lambda_R}$, where by $\lp{t}$ we mean
the ring of Laurent polynomials in one variable $t$ over $\Z$.  Define
$\tau=\alpha_R(t)\in\Lambda_R$. By abuse of notation, we also denote
the greater (in absolute value) of the two eigenvalues of $R$ by
$\tau$. Any $p\in\lp{t}$ can be expressed uniquely in the form
$\sum_{i=m}^Mp_it^i$ where $m$, $M$ and all the $p_i\in\Z$ and neither
$p_m$ nor $p_M$ is $0$. Denote by $m(p)$ and $M(p)$ the values of $m$
and $M$ in this expression and define the \emph{length} of $p$,
$\Abs{p}$, to be $\sum_{i=m}^M\abs{p_i}$. For $\lambda\in\Lambda_R$,
define $\overline{\lambda}$ to be the image of $\lambda$ under the
unique nontrivial (ring) involution of $\Lambda_R$. Correspondingly,
for $p\in\lp{t}$, define $\overline{p}$ to be the image of $p$ under
the corresponding involution of $\lp{t}$, which takes $t$ to $\pm1/t$,
the sign depending on the determinant of $R$. This choice of sign
serves to guarantee that
$\alpha_R(\overline{p})=\overline{\alpha_R(p)}$.

The action by $R$ gives rise to the structure of a $\lp{t}$-module on
$\Z^2=\spanst{x,y}{[x,y]}$, where $t$ acts by $R$. (We will continue
to denote these generators of $\Z^2$ by $x$ and $y$.) This induces a
$\Lambda_R$-action. Let $M_R$ refer to $\Z^2$ with this
$\Lambda_R$-action.

The map $R$, being hyperbolic, has two real eigenvalues, $\tau$ and
$\pm1/\tau$, corresponding to two eigenvectors $\mathbf{v}_e(R)$ and
$\mathbf{v}_c(R)$, taken to be expanding and contracting
respectively. (Clearly, each of these eigenvectors is defined only up
to a multiplicative constant.) The line spanned by $\mathbf{v}_c(R)$
will be called the \emph{contracting line} and the line spanned by
$\mathbf{v}_e(R)$ the \emph{expanding line}. For any $z\in M_R$, let
$d_c(z)$ be the (Euclidean) distance of $z$ from the contracting
line. Then the following diagram commutes:
\[
\begin{CD}M_R@>R\cdot>>M_R\\@Vd_cVV@Vd_cVV\\\R_{\ge0}@>\cdot\abs{\tau}>>\R_{\ge0}.\end{CD}
\]
Similarly, let $d_e(z)$ be the distance of $z$ from the expanding line, so that the following diagram commutes:
\[
\begin{CD}M_R@>R\cdot>>M_R\\@Vd_eVV@Vd_eVV\\\R_{\ge0}@>\div\abs{\tau}>>\R_{\ge0}.\end{CD}
\]

\begin{lem}\label{lower}
Let $\Lambda_R=\Z[t]/p_R(t)$, where $R$ is a hyperbolic automorphism
of $\Z^2$. Then there exist $C_1$, $C_2>0$ with the following
property. Let $\lambda\in\Lambda_R$ and let $\alpha_R$ represent as
usual the natural projection of $\lp{t}$ to $\Lambda_R$. Let $p$,
$q\in{\alpha_R}\inv(\lambda)$, with $p$ a minimal-length element of
this set. Then, for any $z\in M_R$,
\[
d_c(\lambda z)<\abs{\tau}^{M(q)}(C_1(\Abs{q}-\Abs{p})+C_2)d_c(z)
\]
and
\[
d_c(\overline{\lambda}z)<\abs{\tau}^{-m(q)}(C_1(\Abs{q}-\Abs{p})+C_2)d_c(z).
\]
Furthermore,
\[
d_e(\lambda z)<\abs{\tau}^{-m(q)}(C_1(\Abs{q}-\Abs{p})+C_2)d_e(z)
\]
and
\[
d_e(\overline{\lambda}z)<\abs{\tau}^{M(q)}(C_1(\Abs{q}-\Abs{p})+C_2)d_e(z)
\]
\end{lem}

\begin{proof}
We prove only the first statement; the second follows by applying the
involution $\overline{\phantom{p_1}}$ throughout the proof and the
third and fourth then follow by replacing $R$ with $R\inv$ (note that
$\abs{\tr R\inv}=\abs{\tr R}$ since $\abs{\det R}=1$).

Letting brackets represent the greatest integer function, we have $D$
depending only on $R$ such that
\[
\sum_{i=-\infty}^\infty[Dq_i]\le\Abs{q}-\Abs{p},
\]
for otherwise we could apply some multiple of the equation
$t^2\pm1=(\tr R)t$ to $q$ more than $\Abs{q}-\Abs{p}$ times,
shortening $q$ each time without changing $\alpha_R(q)$. This would be
a contradiction since $p$ is of minimal length in
${\alpha_R}\inv(\lambda)$. If $\abs{\tr R}=1$ or $2$ (in which case
the constant term has negative sign, since $R$ is hyperbolic), we must
instead use the equations $3=t^2+t^{-2}$ and $6=t^2+t^{-2}$
respectively. Thus, for any hyperbolic $R$,
\begin{multline*}
\frac{d_c(\lambda
z)}{d_c(z)}\\\le\frac{\Abs{q}-\Abs{p}+1}{D}\tau^{M(q)}+\sum_{i=0}^\infty\frac{1}{D}\abs{\tau}^{M(q)-i-1}\\=\frac{\Abs{q}-\Abs{p}}{D}\abs{\tau}^{M(q)}+\frac{1}{D}\sum_{i=0}^\infty\abs{\tau}^{M(q)-i}\\=\abs{\tau}^{M(q)}\left(\frac{\Abs{q}-\Abs{p}}{D}+\frac{1}{D}\frac{1}{1-\abs{\tau}\inv}\right),
\end{multline*}
as claimed.
\end{proof}

If $p_1$, $p_2\in\lp{t}$, define the \emph{length} of $(p_1,p_2)$ to be
$\Abs{p_1}+\Abs{p_2}$. We mean by $M(p_1,p_2)$ the greater of $M(p_1)$ and
$M(p_2)$ and by $m(p_1,p_2)$ the lesser of $m(p_1)$ and $m(p_2)$.

Observe that conjugation by $c$, denoted by $\rho$, sends $a^ib^j$ to
$a^{i'}b^{j'}$, where $(i',j')=R(i,j)$. Let
$\func[\leftrightarrow]{\beta}{G_R'}{M_R}$ send $a$ to $x$ and $b$ to
$y$. Then $\beta(a^c)=R(\beta(a))$ and $\beta(b^c)=R(\beta(b))$, so
this diagram commutes:
\[
\begin{CD}G_R'@>\cdot^c>>G_R'\\@V\beta VV@V\beta VV\\M_R@>R\cdot>>M_R.\end{CD}
\]
Furthermore, $\beta$ induces an obvious one-to-one length-preserving
correspondence $\theta$ between $N\ab$ and $\lp{t}^2$ such that this
diagram commutes:
\[
\begin{CD}N\ab@>\theta>>\lp{t}^2\\@V\pi VV@VV(p_1,p_2)\mapsto p_1x+p_2yV\\G_R'@>\beta>>M_R.\end{CD}
\]

\begin{defin}
If $v=\sum_{j=1}^nu_{k_j}^{c^{i_j}}$ is a word in $\Z^T$, let
$M(v)$ and $m(v)$ be the maximum and minimum of the $i_j$.
\end{defin}

Proposition~\ref{bdiff} will follow from Proposition~\ref{ll} and from
the following

\begin{prop}\label{taubd}
Let $G_R=\Z^2\rtimes_R\span{c}$, where $R$ is a hyperbolic
automorphism of $\Z^2$. Then there are $D_2>1$ and
$D_1>D_2\ln\abs{\tau}/4$ with the following property. Let $g=uc^z\in
G_R$, where $u\in G_R'$. Let $v$, $v'\in N\ab$ be elements of
$\pi\inv(u)$, with $v$ of minimal length. (We remind the reader that
$\pi$ is the projection from $N\ab$ to $G_R'$.)  Then
$\abs{\tau}^{2M(v)-2M(v')}<D_1(l(v')-l(v))+D_2$ and
$\abs{\tau}^{2m(v')-2m(v)}<D_1(l(v')-l(v))+D_2$.
\end{prop}

The proof follows easily from the following proposition about $M_R$
(whose proof we postpone) and from Lemma~\ref{lower}.

\begin{prop}\label{upper}
Let $M_R$ denote $\Z^2$ equipped with an action by $\lp{t}$, where $t$
acts by some hyperbolic automorphism $R$ of $\Z^2$. Let $x$ and $y$
denote the standard generators of $M_R$ as an abelian group. Then
there exists $e>0$ with the following property. Let $z\in M_R$ and
$p_1$, $p_2\in\lp{t}$ such that $p_1x+p_2y=z$ and
$\Abs{p_1}+\Abs{p_2}$ is minimal subject to this condition. Then
$d_c(z)=d_c(p_1x+p_2y)>e\abs{\tau}^{M(p_1,p_2)}$ and
$d_c(\overline{p_1}x+\overline{p_2}y)>e\abs{\tau}^{-m(p_1,p_2)}$.
\end{prop}

\begin{proof}[Proof of Proposition~\ref{taubd}]
Let $(p_1,p_2)=\theta(v)$ and let $z=\beta(u)=p_1x+p_2y$. Then
$(p_1,p_2)$ is in fact of minimal length such that $p_1x+p_2y=z$
(since $\theta$ is length-preserving and $v$ is a minimal-length
element of $\pi\inv(u)$), so by Proposition~\ref{upper}
$e\abs{\tau}^{M(p_1,p_2)}<d_c(z)$ and
$e\abs{\tau}^{-m(p_1,p_2)}<d_c(\overline{p_1}x+\overline(p_2)y)$.

If similarly $(p_1',p_2')=\theta(v')$ then
$p_1'x+p_2'y=\beta(\pi(v'))=\beta(u)=z$. Then by Lemma~\ref{lower} we
have
\begin{multline*}
e\abs{\tau}^{M(p_1,p_2)}\\<d_c(z)\le
d_c(p_1x)+d_c(p_2y)\\<\abs{\tau}^{M(p_1',p_2')}(C_1(\Abs{p_1'}+\Abs{p_2'}-\Abs{p_1}-\Abs{p_2})+2C_2)\\\cdot\max(d_c(x),d_c(y))
\end{multline*}
and
\begin{multline*}
e\abs{\tau}^{-m(p_1,p_2)}\\<d_c(\overline{p_1}x+\overline(p_2)y)\le
d_c(\overline{p_1}x)+d_c(\overline{p_2}y)\\<\abs{\tau}^{-m(p_1',p_2')}(C_1(\Abs{p_1'}+\Abs{p_2'}-\Abs{p_1}-\Abs{p_2})+2C_2)\\\cdot\max(d_c(x),d_c(y)),
\end{multline*}
where $C_1$ and $C_2$ depend only on $R$ as in Lemma~\ref{lower}.
Thus
\begin{multline*}
\abs{\tau}^{M(p_1,p_2)-M(p_1',p_2')}\\<(C_1(\Abs{p_1'}+\Abs{p_2'}-\Abs{p_1}-\Abs{p_2})+2C_2)\max(d_c(x),d_c(y))
\end{multline*}
and
\begin{multline*}
\abs{\tau}^{m(p_1',p_2')-m(p_1,p_2)}\\<(C_1(\Abs{p_1'}+\Abs{p_2'}-\Abs{p_1}-\Abs{p_2})+2C_2)\max(d_c(x),d_c(y)).
\end{multline*}

Clearly, $M(v)=M(p_1,p_2)$, $m(v)=m(p_1,p_2)$ and
$l(v)=\Abs{p_1}+\Abs{p_2}$, and similarly for $v'$ and
$(p_1',p_2')$. This gives us that
\begin{multline*}
\abs{\tau}^{M(v)-M(v')}=\abs{\tau}^{M(p_1,p_2)-M(p_1',p_2')}\\<(C_1(\Abs{p_1'}+\Abs{p_2'}-\Abs{p_1}-\Abs{p_2})+2C_2)\max(d_c(x),d_c(y))\\=(C_1(l(v')-l(v))+2C_2)\max(d_c(x),d_c(y))
\end{multline*}
and similarly
\[
\abs{\tau}^{m(v')-m(v)}<(C_1(l(v')-l(v))+2C_2)\max(d_c(x),d_c(y)).
\]
(If the $D_2$ resulting from the above expressions is not greater than
$1$, we may clearly increase it without destroying the result. We may
similarly increase $D_1$ to be greater than $D_2\ln\abs{\tau}/4$.)
\end{proof}

\begin{proof}[Proof of Proposition~\ref{bdiff}]
Let $w\in F_S$ be a minimal-length element of
$\phi\inv(v)\cap\sigma_R\inv(g)$ and let $w'\in F_S$ be a
minimal-length element of $\sigma_R\inv(g)$. By definition,
$\abs{g}=l(w')$.

The bound $l(w')\le l(w)$ is obvious, so we need only prove $l(w)-D\le
l(w')$. Let $v'=\phi(w')$. We note that $\pi(v')=u=\pi(v)$.  It is
clear from the definition of $M$ and $m$ that $M(v)\ge m(v)$, so we
cannot have both $M(v)$ negative and $m(v)$ positive, and similarly
for $M(v')$ and $m(v')$. We may thus assume without loss of generality
that $m(v)$ and $m(v')$ are both nonpositive. we will further assume
that $M(v)$ and $M(v')$ are both nonnegative; the other cases are
easier and are left to the reader. In this case, we have by
Lemma~\ref{ll}
\begin{multline*}
[l(w)-l(w')]-[l(v)-l(v')]\\\le2[M(v)-M(v')]+2[m(v')-m(v)]\\{}+\min(\abs{z-M(v)}-M(v),\abs{z-m(v)}+m(v))\\{}-\min(\abs{z-M(v')}-M(v'),\abs{z-m(v')}+m(v')).
\end{multline*}
Noting that $\abs{z-M(v)}-M(v)\ge\abs{z-m(v)}+m(v)$ if and only if
$\abs{z-M(v')}-M(v')\ge\abs{z-m(v')}+m(v')$ if and only if $z\ge0$, we
find using the triangle inequality that
\begin{multline*}
[l(w)-l(w')]-[l(v)-l(v')]-2[M(v)-M(v')]-2[m(v')-m(v)]\\\le\max(\abs{z-M(v)}-\abs{z-M(v')}+M(v')-M(v),\\\abs{z-m(v)}-\abs{z-m(v')}+m(v)-m(v'))\\\le\max(2\abs{M(v')-M(v)},2\abs{m(v)-m(v')})\\\le2\max(M(v)-M(v'),0)+2\max(m(v')-m(v),0).
\end{multline*}
By Proposition~\ref{taubd} and using that $l(v')\ge l(v)$ since $v$ is
given to be of minimal length in $\pi\inv(u)$, we get
\begin{multline*}
[l(w)-l(w')]-[l(v)-l(v')]\\\le4\max(M(v)-M(v'),0)+4\max(m(v')-m(v),0)\\<2\max(\log_{\abs{\tau}}[D_1(l(v')-l(v))+D_2],0)+2\max(\log_{\abs{\tau}}[D_1(l(v')-l(v))+D_2],0)\\=4\log_{\abs{\tau}}[D_1(l(v')-l(v))+D_2].
\end{multline*}
Rearranging yields that
$l(w)-l(w')\le4\log_{\abs{\tau}}[D_1(l(v')-l(v))+D_2]-(l(v')-l(v))$.

It now remains to show that $f(x)=4\log_{\abs{\tau}}(D_1x+D_2)-x$
attains an absolute maximum value on the domain $[0,\infty)$, which
can be done by differential calculus. (This absolute maximum value
will depend only on $R$ since $\tau$, $D_1$ and $D_2$ depend only on
$R$.)  We compute $f'(x)=4D_1/[(D_1x+D_2)\ln\abs{\tau}]-1$ and observe
that it is a strictly decreasing function on the domain of interest,
so (by the first derivative test for absolute maxima), if it has a
zero, the absolute maximum of $f$ is attained there. We find that
$\lim_{x\to\infty}f'(x)=-1$ while $f'(0)=4D_1/(D_2\ln\abs{\tau})-1>0$
since we took $D_1>D_2\ln\abs{\tau}/4$. Since $f'$ is continuous on
$[0,\infty)$, we are done.
\end{proof}

The proof of Proposition~\ref{bdiff} thus reduces to that of
Proposition~\ref{upper}. This in turn follows from the following two
results.

\begin{lem}\label{upperres}
Let $M_R$ denote $\Z^2$ equipped with an action by $\lp{t}$, where $t$
acts by some hyperbolic automorphism $R$ of $\Z^2$. Let $x$ and $y$
denote the standard generators of $M_R$ as an abelian group. Let
$\Lambda_R$ denote the quotient of $\lp{t}$ by the charcteristic
polynomial of $R$. Let $\lambda_1$, $\lambda_2\in\Lambda_R$ and $p_1$,
$p_2\in\lp{t}$ be minimal-length elements of $\alpha_R\inv(\lambda_1)$
and $\alpha_R\inv(\lambda_2)$ respectively. Then
\[
6(d_c(x)+d_c(y))\abs{\tau}^{M(p_1,p_2)+1}>d_c(p_1x+p_2y)
\]
and
\[
6(d_c(x)+d_c(y))\abs{\tau}^{-m(p_1,p_2)+1}>d_c(\overline{p_1}x+\overline{p_2}y).
\]
\end{lem}

\begin{proof}
We prove the first statement of the lemma; the second will follow by
applying the involution $\overline{\phantom{p_1}}$ throughout the
proof.

We know, assuming $\abs{\tr R}>2$, that no coefficient in $p_1$ or
$p_2$ can have absolute value greater than $\abs{\tr R}+2$, since
otherwise we could apply the equation $\pm t^2-(\tr R)t\pm1=0$ to
produce a shorter element of $\alpha_R\inv(\lambda_1)$ or
$\alpha_R\inv(\lambda_2)$. If $\abs{\tr R}=1$ or $2$ the same result
holds, since instead we use the equation $t^2+t^{-2}=3$ or $6$,
respectively. Thus, in either case, given $z\in M_R$, both
$d_c(p_1z)/d_c(z)$ and $d_c(p_2z)/d_c(z)$ are
\[
\le\sum_{i=-\infty}^{M(p_1,p_2)}(\abs{\tr R}+2)\abs{\tau}^i=(\abs{\tr R}+2)\frac{\abs{\tau}^{M(p_1,p_2)+1}}{\abs{\tau}-1}.
\]
If we denote $\abs{\tr R}$ by $r$ then this equals
\[
\frac{2(r+2)\abs{\tau}^{M(p_1,p_2)+1}}{r+\sqrt{r^2\mp4}-2}<6\abs{\tau}^{M(p_1,p_2)+1}
\]
since either $r>2$ or the sign inside the square root is positive. It
follows that
\[
d_c(p_1x+p_2y)\le
d_c(p_1x)+d_c(p_2y)<6(d_c(x)+d_c(y))\abs{\tau}^{M(p_1,p_2)+1},
\]
as claimed.
\end{proof}

\begin{prop}\label{upperhead}
Let $M_R$ denote $\Z^2$ equipped with an action by $\lp{t}$, where $t$
acts by some hyperbolic automorphism $R$ of $\Z^2$. Let $x$ and $y$
denote the standard generators of $M_R$ as an abelian group and let
$P$ and $Q$ be positive numbers. Then there exists a natural number
$n$ with the following property. Let $z\in M_R$ and $(p_1,p_2)$ be a
minimal-length element of $\lp{t}^2$ such that $p_1x+p_2y=z$. Then
$(Pd_c(x)+Qd_c(y))\abs{\tau}^{M(p_1,p_2)-n}<d_c(z)$.
\end{prop}

We postpone the proof and instead first give the

\begin{proof}[Proof of Proposition~\ref{upper}]
As in the proof of Lemma~\ref{upperres}, we prove the first statement
of the lemma; the second will follow by applying the involution
$\overline{\phantom{p_1}}$ throughout the proof.

If $p\in\lp{t}$, we let $\head_m(p)$ denote the sum of those terms of
$p$ of degree at least $m$, and $\tail_m(p)$ the sum of all remaining
terms. (Thus $p=\head_m(p)+\tail_m(p)$.) If $m$ is taken to be
$M(p_1,p_2)-n$, with $n$ given by Proposition~\ref{upperhead} with
$P=7$ and $Q=6$, we have
$d_c(\head_m(p_1)x+\head_m(p_2)y)>(7d_c(x)+6d_c(y))\abs{\tau}^{M(p_1,p_2)-n}$
by that proposition. However,
$M(\tail_m(p_1),\tail_m(p_2))=m-1$. Thus, by Lemma~\ref{upperres}, we
have
\begin{multline*}
d_c(\tail_m(p_1)x+\tail_m(p_2))<6(d_c(x)+d_c(y))\abs{\tau}^m\\=6(d_c(x)+d_c(y))\abs{\tau}^{M(p_1,p_2)-n},
\end{multline*}
since certainly if $(p_1,p_2)$ is of minimal length such that
$p_1x+p_2y=z$ then each of the $p_i$ is of minimal length in
$\alpha_R\inv(\alpha_R(p_i))$. (We again remind the reader that
$\alpha_R$ is the projection from $\lp{t}$to $M_R$.) It follows that
\begin{multline*}
d_c(z)\ge
d_c(\head_m(p_1)x+\head_m(p_2)y)-d_c(\tail_m(p_1)x+\tail_m(p_2)y)\\>d_c(x)\abs{\tau}^{M(p_1,p_2)-n}.
\end{multline*}
We may thus take $e$ to be $d_c(x)\abs{\tau}^{-n}$.
\end{proof}

It remains now (for the proof of Proposition~\ref{bdiff}) to prove
Proposition~\ref{upperhead}. It will follow from the following three
results.

\begin{lem}\label{lin}
Let $M_R$ denote $\Z^2$ equipped with an action by $\lp{t}$, where $t$
acts by some hyperbolic automorphism $R$ of $\Z^2$. Let $x$ and $y$
denote the standard generators of $M_R$ as an abelian group and let
$d\in\R$. Then there exist $C_1>0$ and $C_2\in\R$ with the following
property. Let $p_1$, $p_2\in\lp{t}$ be such that $d_c(p_1x+p_2y)\le d$
and the first coordinate of $p_1x+p_2y\in M_R$ is $w$. Assume that
$m(p_1,p_2)\ge0$. Then $\Abs{p_1}+\Abs{p_2}\ge C_1\abs{w}+C_2$.
\end{lem}

\begin{proof}
Since $m(p_1,p_2)\ge0$, $d_e(p_1x+p_2y)\le\Abs{p_1}+\Abs{p_2}$, since
$p_1x+p_2y$ is the sum of $\Abs{p_1}+\Abs{p_2}$ vectors $v_i$ ($=t^jx$
or $t^jy$, $j\in\Z$) such that $d_e(v_i)\le1$. However, the nearest
point to $p_1x+p_2y$ on the contracting line has first coordinate of
absolute value at least $\abs{w}-d$. This point is thus at distance at
least $C_1(\abs{w}-d)$ from the expanding line, where $C_1$ is s
positive constant depending only on $R$. Since $p_1x+p_2y$ is within
distance $d$ of this point, $d_e(p_1x+p_2y)\ge C_1(\abs{w}-d)-d$, as
claimed. The result follows by transitivity if we let $C_2=-d(C_1+1)$.
\end{proof}

We postpone the proof of the next two results.

\begin{prop}\label{log}
Let $M_R$ denote $\Z^2$ equipped with an action by $\lp{t}$, where $t$
acts by some hyperbolic automorphism $R$ of $\Z^2$. Let $x$ and $y$
denote the standard generators of $M_R$ as an abelian group and let
$d\ge0$. Then there exist $D_1$, $D_2$ and $D_3$ such that the
following holds. Let $z\in M_R$ be such that $d_c(z)\le d$ and the
first coordinate of $z$ is $w$. Then there exist $p_1$, $p_2\in\lp{t}$
such that $p_1x+p_2y=z$ and $\Abs{p_1}+\Abs{p_2}\le
D_1+\max(0,D_2\ln(D_3\abs{w}))$.
\end{prop}

\begin{prop}\label{finite}
Let $M_R$ denote $\Z^2$ equipped with an action by $\lp{t}$, where $t$
acts by some hyperbolic automorphism $R$ of $\Z^2$. Let $x$ and $y$
denote the standard generators of $M_R$ as an abelian group and let
$z\in M_R$. Let $P$ be the set of all $(p_1,p_2)\in\lp{t}^2$ such that
$p_1x+p_2y=z$. Then $P$ contains only finitely many elements of
minimal length.
\end{prop}

\begin{proof}[Proof of Proposition~\ref{upperhead}]
It will suffice to show for some $n$ that, if $M(p_1,p_2)\ge n$,
then $Pd_c(x)+Qd_c(y)<d_c(z)$.

Suppose the conclusion does not hold for $z$, $p_1$ and $p_2$; that
is, suppose that $Pd_c(x)+Qd_c(y)\ge d_c(z)$ where $z=p_1x+p_2y$. Then
applying Lemma~\ref{lin} with $d=Pd_c(x)+Qd_c(y)$ (which, it should be
noted, depends only on $R$, $P$ and $Q$) yields $C_1>0$ and $C_2$ such
that
\[
\Abs{p_1}+\Abs{p_2}\ge C_1\abs{w}+C_2,
\]
where $w$ is the first coordinate of $z$.  However, since $p_1$ and
$p_2$ were taken so that $\Abs{p_1}+\Abs{p_2}$ would be minimal
subject to the condition $p_1x+p_2y=z$, applying Proposition~\ref{log}
with the same value of $d$ gives $D_1$, $D_2$ and $D_3$ such that
$\Abs{p_1}+\Abs{p_2}\le D_1+\max(0,D_2\ln(D_3\abs{w}))$. Thus
\[
D_1+\max(0,D_2\ln(D_3\abs{w}))\ge C_1\abs{w}+C_2.
\]
Since the left-hand side depends logarithmically on $\abs{w}$ while
the right-hand side is increasing and depends linearly on it, there
exists $w'$ dependent only on $R$, $P$ and $Q$ (via $D_1$, $D_2$,
$D_3$, $C_1$ and $C_2$) such that $w<w'$. Thus $z$ is a lattice point
constrained to lie within a compact set, since it is confined to the
intersection of finite-width bands about the $y$-axis and the
contracting line. There are only finitely namy such points. Thus there
are also only finitely many possibilities for $p_1$ and $p_2$, by
Proposition~\ref{finite}. The result follows by taking $n$ to be $1$
more than the maximal $M(p_1,p_2)$ for any of these finitely many
possibilities.
\end{proof}

Propsition~\ref{log} will follow from the following

\begin{lem}\label{powers}
Let $R$ be a hyperbolic automorphism of $\Z^2$. Let $S$ be the set of
all integers appearing in the first row of the matrix for some
nonnegative power of $R$. Then there exist $C_1$, $C_2$ and $C_3$ such
that any $n\in\Z$ can be expressed as a word in the elements of $S$ of
length at most $C_1+\max(0,C_2\ln(C_3\abs{n})$.
\end{lem}

\begin{proof}
The matrix for $R$ consists of two column vectors, which we denote
$(p,q)$ and $(r,s)$. Each can be decomposed as the sum of a (real)
vector in the expanding line and one in the contracting line; thus
$(p,q)=(p_e,q_e)+(p_c,q_c)$ and $(r,s)=(r_e,s_e)+(r_c,s_c)$ where
$(p_e,q_e)$ and $(r_e,s_e)$ lie in the expanding line and $(p_c,q_c)$
and $(r_c,s_c)$ in the contracting line. If we denote the expanding
eigenvalue of $R$ by $\tau$ then we find that
\[
S=\setst{p_m=p_e\tau^m\pm\frac{p_z}{\tau^m},r_m=r_e\tau^m\pm\frac{r_z}{\tau^m}}{m\ge0},
\]
where the $\pm$ signs depend on the determinant of $R$ and the parity
of $m$.

Let $C_1$ be the maximal length of any minimal-length word in the
elements of $S$ representing a number of absolute value
$\le\abs{p_e}$. Let $C_2$ be $\abs{\tau}(1+\abs{p_z})/\ln\abs{tau}$
and $C_3$ be $\abs{\tau/p_e}$.

Consider $n,m\in\Z$ with $\abs{n}\le\abs{p_e\tau^m},m>0$. Then there
must exist $k\in\Z$ with $\abs{k}<\abs{\tau}$ such that
$\abs{n+kp_e\tau^{m-1}}\le\abs{p_e\tau^{m-1}}$. It follows that
\[
\abs{n+kp_{m-1}}\le\abs{n+kp_e\tau^{m-1}}+\abs{\frac{kp_c}{\tau^{m-1}}}\le\abs{p_e\tau^{m-1}}+\abs{kp_c}\le\abs{p_e\tau^{m-1}}+\abs{kp_c},
\]
so there exists $l\in\Z$ with $\abs{l}\le\abs{kp_c}<\abs{\tau p_c}$
such that $\abs{n+kp_{m-1}+lp_0}\le\abs{p_e\tau^{m-1}}$. (Note that
$p_0=1$.) The result follows by induction on $m$.
\end{proof}

\begin{proof}[Proof of Proposition~\ref{log}]
By Lemma~\ref{powers} applied to the hyperbolic automorphism $R\inv$,
there exist $C_1$, $C_2$ and $C_3$ depending only on $R$ and $p_1$,
$p_2'\in\lp{t}$ such that $p_1x+p_2'y$ has first coordinate $w$ and
\[
\Abs{p_1}+\Abs{p_2'}\le C_1+\max(0,C_2\ln(C_3\abs{w})).
\]
Furthermore, $M(p_1,p_2')\le0$. Then
\[
d_c(p_1x+p_2'y)\le C_1+\max(0,C_2\ln(C_3\abs{w}))
\]
since $p_1x+p_2'y$ is the sum of $C_1+\max(0,C_2\ln(C_3\abs{w}))$
vectors $v_i$ ($=t^jx$ or $t^jy$, $j\in\Z$) such that
$d_c(v_i)\le1$. Thus $p_1x+p_2'y$ is within
$\sqrt{2}[C_1+\max(0,C_2\ln(C_3\abs{w}))]$ vertically of the point on
the contracting line with first coordinate $w$, so within
$\sqrt{2}[C_1+\max(0,C_2\ln(C_3\abs{w}))+d]$ of $z$. (It is here that
we use the condition that $\abs{zeta}\le1$.)  Defining $p_2$ to be
\[
p_2'+(\text{some appropriate constant term}),
\]
we see that $D_1=(1+\sqrt{2})C_1+\sqrt{2}d$, $D_2=(1+\sqrt{2})C_2$ and
$D_3=C_3$ give what we seek.
\end{proof}

The proof of Proposition~\ref{bdiff} has now been reduced to that of
Proposition~\ref{finite}, which follows from the following

\begin{lem}\label{gaps}
Let $M_R$ denote $\Z^2$ equipped with an action by $\lp{t}$, where $t$
acts by some hyperbolic automorphism $R$ of $\Z^2$. Let $x$ and $y$
denote the standard generators of $M_R$ as an abelian group and let
$z\ne\mathbf{0}\in M_R$ and $l\in\N$.  Then there exists an integer
$N=N(z,l,R)$ with the following property. Let $p_1$, $p_2\in\lp{t}$
such that
\[
p_1x+p_2y=z
\]
and $\Abs{p_1}+\Abs{p_2}=l$. Then at least one of $p_1$ and $p_2$ must
have a nonzero term whose degree has absolute value less than $N$.
\end{lem}

\begin{proof}
Express $z$ (viewed as an element of $\Z^2\subset\R^2$) as
$z=z_e\mathbf{v}_e(R)+z_c\mathbf{v}_c(R)$, $z_e$, $z_c\in\R$. (We
remind the reader that $\mathbf{c}_e(R)$ and $\mathbf{v}_c(R)$
generate the expanding and contracting eigenspaces respectively.)
Since both eigenvalues of $R$ are irrational, no nonzero real multiple
of $\mathbf{v}_c(R)$ lies in the $\Z^2$ lattice. In particular, since
$z\ne\mathbf{0}$, $z_c\mathbf{v}_c(R)\notin\Z^2$. Let $D>0$ be the
(Euclidean Hausdorff) distance between $z_c\mathbf{v}_c(R)$ and
$\Z^2$, and let $N>\log_{\abs{\tau}}(2l/D)$.

If $p\in\lp{t}$, let $p_{\ge m}\in\lp{t}$ denote the sum of all terms
of $p$ of degree $\ge m$, and analogously for $p_{\le m}$. Suppose
(for a contradiction) $p_1={p_1}_{\ge N}+{p_1}_{\le-N}$ and similarly
for $p_2$. Then
\[
p_1x+p_2y={p_1}_{\ge N}x+{p_2}_{\ge N}y+{p_1}_{\le-N}x+{p_2}_{\le-N}y.
\]
We define ${z_e}_{\ge N}$ and ${z_c}_{\ge N}\in\R$ such that
\[
{p_1}_{\ge N}x+{p_2}_{\ge N}y={z_e}_{\ge N}\mathbf{v}_e(R)+{z_c}_{\ge
N}\mathbf{c}(R),
\]
and similarly for ${z_e}_{\le-N}$ and ${z_c}_{\le-N}$. Then each term
$a_it^ix$ or $a_it^iy$ of ${p_1}_{\ge N}x+{p_2}_{\ge N}y$ contributes
at most $\abs{a_i}/\tau^N$ of Euclidean length to ${z_c}_{\ge
N}\mathbf{v}_c(R)$. Since te sum of all the $a_i$ is at most $l$,
${z_c}_{\ge N}\mathbf{v}_c(R)$ has Euclidean length at most
$l/\tau^N$. Similarly, ${z_e}_{\le-N}\mathbf{v}_e(R)$ has Euclidean
length at most $l/\tau^N$. Then
\[
{p_1}_{\le -N}x+{p_2}_{\le-N}y=z_c\mathbf{v}_c(R)-{z_c}_{\ge
  N}\mathbf{v}_c(R)+{z_e}_{\le-N}\mathbf{v}_c(R)\in\Z^2
\]
lies within Euclidean distance $2l/\tau^N<D$ of
$z_c\mathbf{v}_c(R)$. But this is impossible by the definition of $D$.
\end{proof}

\begin{proof}[Proof of Proposition~\ref{finite}]
The result is trivial for $z=\mathbf{0}$, since then $P=\set{(0,0)}$,
so we assume $z\ne\mathbf{0}$.

Suppose $P$ contains infinitely many minimal-length elements. Let $l$
be the length of these minimal-length elements. Then they each contain
some term in which $t$ is raised to an exponent of absolute value
$<N=N(z,l,R)$, by Lemma~\ref{gaps}. Thus each minimal-length element
of $P$ comes from a length-($l-1$) minimal-length expression for one
of the finitely many possibilities $z\pm t^{-N+1}x$, $z\pm t^{-N+1}y$,
$z\pm t^{-N+2}x$, \dots, $z\pm t^{N-1}y$. The result follows by
induction on $l$.
\end{proof}

This concludes the proof of Proposition~\ref{bdiff}. We now need only
to prove Proposition~\ref{flat}.

\begin{lem}\label{distort}
Let $G_R=\Z^2\rtimes_R\Z$, where $R$ is a hyperbolic automorphism of
$\Z^2$, and let $z=(z_1,z_2)\in G_R'$ be such that $\abs{z_1}$,
$\abs{z_2}<\abs{\tr R}^m$. Then there exists a word
$w\in\sigma_R\inv(\beta\inv(z))$ such that
$l(w)\le2^{m+1}+4m-1$. (Recall that $\beta$ is the correspondence
between $G_R'$ and $\Z^2$ and $\sigma_R$ is the surjection from $F_S$
to $G_R$.) If $\abs{\tr R}=1$ then if $\abs{z_1}$, $\abs{z_2}<3^m$
then there exists $w\in\sigma_R\inv(\beta\inv(z))$ such that
$l(w)\le2^{m+1}+8m-1$. If $\abs{\tr R}=2$ then if $\abs{z_1}$,
$\abs{z_2}<6^m$ then there exists $w\in\sigma_R\inv(\beta\inv(z))$
such that $l(w)\le2^{m+1}+8m-1$.
\end{lem}

\begin{proof}
We begin with the proof of the first claim.

Since $\abs{z_1}$, $\abs{z_2}<\abs{\tr R}^m$, we can express them in
the form
\[
z_1=\sum_{i=0}^mz_{1i}\abs{\tr R}^i,
\]
\[
z_2=\sum_{i=0}^mz_{2i}\abs{\tr R}^i,
\]
where all the $\abs{z_{1i}}$ and $\abs{z_{2i}}\le\abs{\tr R}/2$. But
since, in $\Lambda_R$, $\abs{\tr R}=\pm\tau\pm\tau\inv$, $\abs{\tr
R}^i=(\pm\tau\pm\tau\inv)^i=\sum_{j=0}^i\pm\binom{i}{j}\tau^{2j-i}$. Thus
we can reexpress $z_1$ as $\sum_{k=-m}^mz_{1k}'\tau^k$, where
\[
\sum_{k=-m}^m\abs{z_{1k}'}\le\sum_{i=0}^m\sum_{j=0}^i\binom{i}{j}=\sum_{i=0}^m2^i=2^{m+1}-1,
\]
and similarly for $z_2$. We thus have
$z=(\sum_{k=-m}^mz_{1k}'t^k)x+(\sum_{k=-m}^mz_{2k}'t^k)y$. Thus there
exists $v=\sum_{j=1}^nu_{k_j}^{c^{i_j}}\in\Z^T$ with
$\pi(v)=\beta\inv(z)$ and all $k_j=1$ or $2$ such that all $i_j$
satisfy $-m\le i_j\le m$ and $n\le2^{m+1}-1$. By Lemma~\ref{ll}, this
gives a $w\in\sigma_R\inv(\beta\inv(z))$ such that that
$l(w)\le2^{m+1}+4m-1$, as claimed.

The proofs of the second and third claims are analogous; just replace
$\abs{\tr R}$ by $3$ and $6$ respectively and the equation $\abs{\tr
R}=\pm\tau\pm\tau\inv$ by the equations $3=\tau^2+\tau^{-2}$ and
$6=\tau^2+\tau^{-2}$ respectively.
\end{proof}

\begin{proof}[Proof of Proposition~\ref{flat}]
We first construct $B_{m,n}$. Recall that $\beta$ is the
correspondence between $G_R'$ and $M_R$ and $\theta$ is that between
$N\ab$ and $\lp{t}^2$. Let $u\in G_R'$ and set $v\in N\ab$ to be a
minimal-length element of $\pi\inv(u)$. Then by Lemma~\ref{lower}
(setting $(q_1,q_2)=(p_1,p_2)=\theta(v)$)
\begin{multline*}
d_c(\beta(u))\le
d_c(p_1x)+d_c(p_2y)\\<\abs{\tau}^{M(\theta(v))}(C_1(\Abs{q_1}+\Abs{q_2}-\Abs{p_1}-\Abs{p_2})+2C_2)\max(d_c(x),d_c(y))\\=2C_2\abs{\tau}^{M(\theta(v))}\max(d_c(x),d_c(y))
\end{multline*}
and
\begin{multline*}
d_e(\beta(u))\\<\abs{\tau}^{-m(\theta(v))}(C_1(\Abs{q_1}+\Abs{q_2}-\Abs{p_1}-\Abs{p_2})+2C_2)\\\cdot\max(d_e(x),d_e(y))\\=2C_2\abs{\tau}^{-m(\theta(v))}\max(d_e(x),d_e(y))
\end{multline*}
where $C_2$ depends only on $R$.  Then if we assume
\[
d_c(\beta(u))>2C_2\abs{\tau}^n\max(d_c(x),d_c(y)),
\]
and
\[
d_e(\beta(u))>2C_2\abs{\tau}^n\max(d_e(x),d_e(y)),
\]
this implies that $M(\theta(v))>n$, $m(\theta(v))<-n$, so $M(v)>n$,
$m(v)<-n$. We let $B_{m,n}$ be the set of $u=a^ib^j\in[G_R,G_R]$ such
that $\abs{i}$, $\abs{j}<\abs{\tr R}^m$,
\[
d_c(\beta(a^ib^j))>2C_2\abs{\tau}^n\max(d_c(x),d_c(y))
\]
and
\[
d_e(\beta(a^ib^j))>2C_2\abs{\tau}^n\max(d_e(x),d_e(y)).
\]
If $\abs{\tr R}=1$ then replace $\abs{\tr R}$ in the preceding
sentence by $3$.

We proceed with the proof of the rest of the proposition, first
assuming $\abs{\tr R}>2$.

Let $K$ and $L$ be positive integers such that
\[
(L+2C_2\abs{\tau}^n)\frac{\max(d_c(x),d_c(y),d_e(x),d_e(y))}{\min(d_c(x),d_e(x))}<K<\abs{\tr R}^m-L;
\]
such $K$ and $L$ clearly exist for sufficiently large $m$ and we can
arrange for $L$ to grow proportionally to $\abs{\tr R}^m$ for $R$ and
$n$ fixed.  Let $z=Kx$.  Addition or subtraction of $x$ or $y$ does
not change $d_c(z)$ by more than $\max(d_c(x),d_c(y))$ nor $d_e(z)$ by
more than $\max(d_e(x),d_e(y))$. Thus any element
$z'=(z_1',z_2')\in\Z^2$ at $\set{x,y}$-distance less than $L$ from $z$
satisfies
\[
d_c(z')>Kd_c(x)-L\max(d_c(x),d_c(y))>2C_2\abs{\tau}^n\max(d_c(x),d_c(y)),
\]
$d_e(z')>2C_2\abs{\tau}^n\max(d_e(x),d_e(y))$. But $z_1'$ and $z_2'$
must then also satisfy $\abs{z_1'}$, $\abs{z_2'}<K+L<\abs{\tr
R}^m$. Then there is $w\in\sigma_R\inv(\beta\inv(z))$ with
$l(w)\le2^{m+1}+4m-1$ by Lemma~\ref{distort}, that is to say
$\abs{\beta\inv(z)}\le2^{m+1}+4m-1$. Then, by Proposition~\ref{bdiff},
$\Abs{\beta\inv(z')}\le2^{m+1}+4m+D-1$. Letting
$r=[L/(2^{m+1}+4m+D-1)]-1$, where the brackets represent the greatest
integer function, the result now follows from Proposition~\ref{fuzz}
if we take $f(z')=\Abs{\beta\inv(z')}$. It is straightforward to check
that $\lim_{m\to\infty}r=\infty$, since $\abs{\tr R}>2$.

The proof for $\abs{\tr R}\le2$ is analogous; just replace $\abs{\tr
R}$ by $3$ or $6$ as needed and use the other half of
Lemma~\ref{distort}.
\end{proof}

\end{document}